\documentclass[preprint]{elsarticle}

\usepackage[a4paper,%
left=1in,right=1in,top=1in,bottom=1in,%
footskip=.25in]{geometry}

\usepackage{amsfonts,amssymb}
\usepackage{graphicx,epstopdf}
\usepackage{tikz}
\usetikzlibrary{automata}
\usepackage{color}
\usepackage{threeparttable}
\usepackage{amsmath}
\usepackage{amsfonts,amssymb}
\usepackage{slashbox} 
\usepackage{epstopdf}
\usepackage{algorithm}
\usepackage[noend]{algpseudocode}
\usepackage{ulem}
\usepackage{etoolbox}
\usepackage{palatino}

\patchcmd{\subequations}{}%
{}{}{}

\allowdisplaybreaks[4]

\usepackage{lineno,hyperref}

\journal{IEEE Intelligent Transportation Systems Conference 2019, Auckland, New Zealand}

\begin{document}

\begin{frontmatter}

\title{\LARGE \bf A Distributed Architecture for Real-time Hybrid Traffic Light Control in Urban Transportation Networks\\
	\normalsize -- Presentation at IEEE ITSS Young Professionals Workshop 2019
}
\tnotetext[mytitlenote]{This work was presented during the IEEE ITSS Young Professionals Workshop affiliated with IEEE Intelligent Transportation Systems Conference 2019 at Auckland, New Zealand.}

\author[i2raddr]{Yicheng Zhang\fnref{i2r}}
\address[i2raddr]{Institute for Infocomm Research (I2R), 1 Fusionopolis Way, 21-01 Connexis, Singapore 138632}
\fntext[i2r]{Yicheng Zhang is affiliated with the Institute for Infocomm Research (I2R) at the Agency for Science, Technology and Research (A*STAR), 1 Fusionopolis Way, 21-01 Connexis, Singapore 138632. Email: {\tt\small yzhang088@e.ntu.edu.sg}
}

%
%

\begin{abstract}
A macroscopic model is proposed to depict the traffic dynamics involved in urban traffic systems. The link dynamics are described based on the cell-transmission model and bounded by the link capacities, while the flow dynamics are proposed based on the discharge headways and saturation flow at intersections. To fulfill the requirement of a closed-loop traffic light control strategy, an approach to estimate the branching ratios at intersections is proposed and simulations show that the convergence would be achieved under constant cyclic flow profiles. Furthermore, a system partitioning approach is proposed based congestion level identification, which is achieved via a machine learning method and a hybrid traffic network control strategy is proposed to integrate different traffic light control schemes together. 
\end{abstract}


\end{frontmatter}


\section{Introduction}

Traffic congestions in urban traffic system is a persisting problem. 
The conventional traffic light control is based on fixed-time strategy. 
In modern urban traffic systems, some transportation responsive strategies are widely used such as SCOOT and SCATS.
However, most traffic light scheduling strategies focus on developing well-tuned off-line solutions.
Recently, more advanced traffic light control strategies are introduced to improve the urban mobility.
Strategies like OPAC, PRODYN, and RHODES introduce optimizations into solving the traffic light control problem. However, those optimization algorithms are with exponentially computational complexity and are not suitable for centralized implementations for large-scale traffic networks. 
Distributed algorithms with system partitioning strategies become a possible solution for the computational issue. 
Partitioning methodologies are first introduced by computer engineers for solving problems like data mining and image segmentation. K-means method is one of the most widely used algorithm in computer engineering as well as in transportation engineering. 
A recent work shows the traffic system partitioning based on a macroscopic fundamental diagram, which is to partition the traffic system based on the trade-off between partitioning (normalized cut algorithm) and merging (agglomerative clustering algorithm).
Traditional distributed traffic light control strategies ignore the coordination among the signals at various intersections, which could not achieve a system-level optimization. 

\section{Main Results}
Solving the traffic light control problem in a large-scale network is time consuming. One straightforward idea is to partition the large-scale network into small regions. Three issues are involved in this procedure. Firstly, the traffic network is a connected graph, how to obtain the partitioning scheme and how to evaluate the quality of the obtained partitioning scheme are the fundamental questions to answer. Secondly, solving all the local problems optimally is not equivalent to solving the global problem optimally, how to fill the gap between the ``local" optimality and global optimality is another issue. Thirdly, the computational complexity involved in this problem comes from the aspiration of fully optimizing the traffic network, but do we really need all intersections to be controlled ``smartly"? 
In our previous work \cite{Zhang2015}, a distributed computational structure was proposed to answer the second question. 
In this design, a learning based congestion identification method is proposed for system partitioning and a hybrid traffic light control strategy is proposed to reduce the computational complexity by integrating different strategies.

\begin{figure}[!htbp]
	\centering
	\includegraphics[width = 0.85\linewidth]{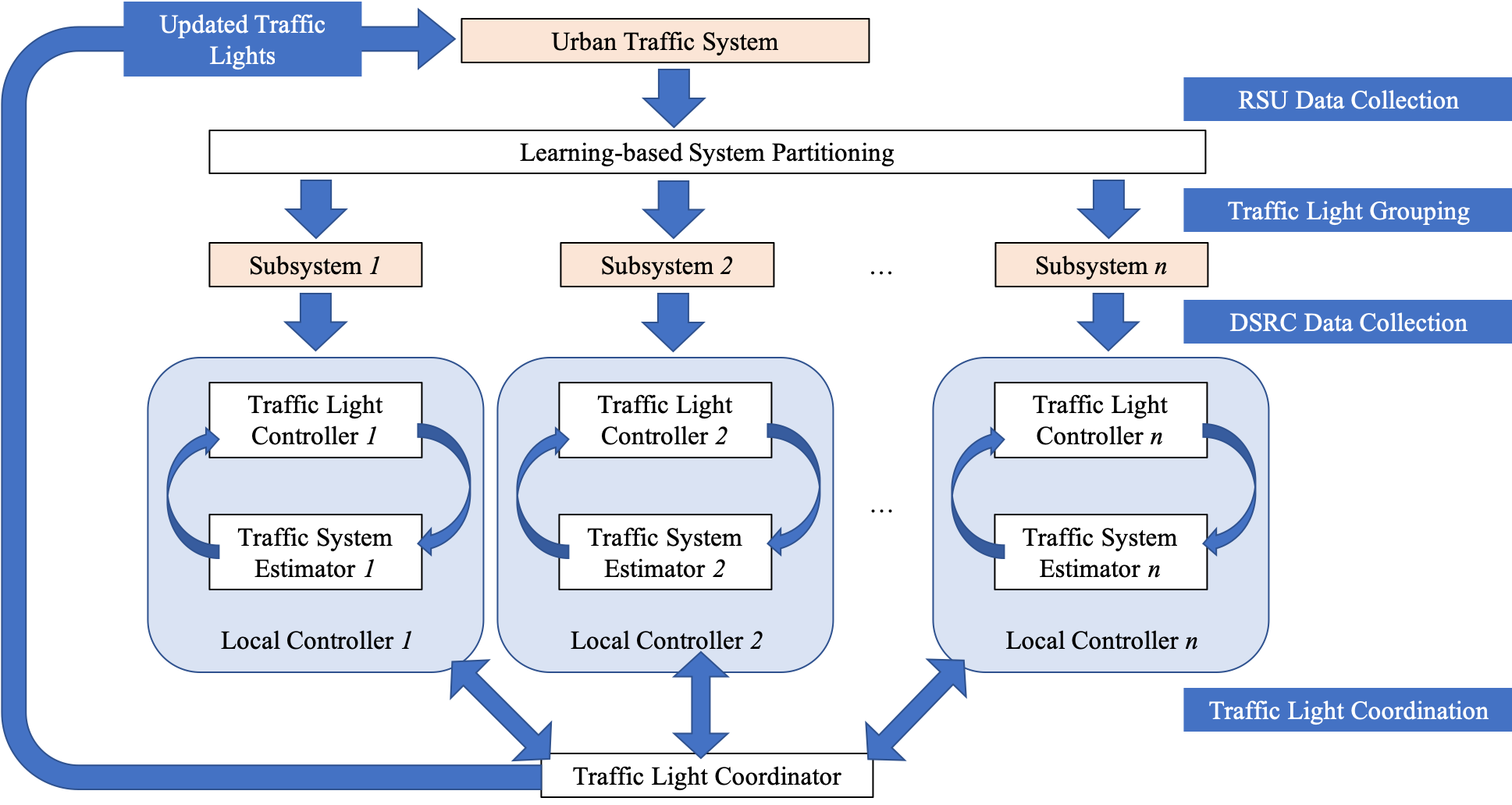}
	\caption{The Flowchart for Hybrid Traffic Light Control Strategy}
	\label{fig_flowchart}
\end{figure}

The flowchart for the hybrid traffic light control strategy is shown in Fig. \ref{fig_flowchart}.

The contributions of this work are summarized as follows: 
(1) A macroscopic urban traffic system model with revised piecewise linear function based flow dynamics and its validation;
(2) A closed-loop traffic light control strategy design and branching ratio estimations based on moving-average model;
(3) Congestion-level based traffic system partitioning;
(4) Design of a distributed traffic management system with hybrid traffic light control strategies.


\subsection{A Learning Based Congestion Identification Method and Traffic System Partitioning}

%
In this work, the traffic flow based method is adopted.
Two 1-to-n recursive neural networks (RNNs) are adopted to predict the congestion levels, i.e., one model to predict the vehicle velocity on each lane and the other model to predict the traffic density. The hyper parameters will be well-trained via historical data. 
The main idea is to learn the congestion levels based on traffic measurements and use the congestion levels as criteria to partition the traffic network. 
The process could be partitioned into three stages, the link congestion level prediction, identification and link clustering. 
The system partitioning information will be used in the design of a hybrid traffic light control scheme, which is shown in the following section. 


\subsection{A Hybrid Traffic Light Control Scheme Based on System Partitioning}
Based on congestion level identification for each link and the link clustering algorithm, the system could be partitioned into subsystems.
Computations for all the subsystems could be executed in parallel. 
Different sub-systems could adopt variant traffic light control strategies based on the traffic requirements. 

Then the question goes to how the consistency of boundary links among different subsystems could be obtained. 
Two scenarios are involved in this problem, i.e., (1) the traffic light control strategies adopted in the two adjacent subsystems are traffic responsive ones, or (2) the traffic light control strategies adopted in the two adjacent subsystems are traffic responsive one and pre-timed scheme, respectively. 
For the second case, the traffic dynamics in the subsystems which are adopted pre-timed scheme are calculated in the first stage. Thus, the estimated outputs from these subsystems are treated as the inputs for the other subsystems which are using the traffic responsive schemes. 
For the first case, similar to what we proposed in \cite{Zhang2015}, a distributed traffic light control scheme based on Lagrangian multiplier method is adopted to achieve the consistency. 


\section{Simulation Results}
\label{sec_simulation}

The proposed distributed hybrid traffic light scheme is simulated based on West Singapore Transportation System, which is built in PTV VISSIM. 
The optimized traffic signal controller is utilized when the congestion level is high. To simplify the integration, a pre-timed scheme is selected for the specific region.

\begin{figure}[!htbp]
	\centering
	\includegraphics[width=\linewidth]{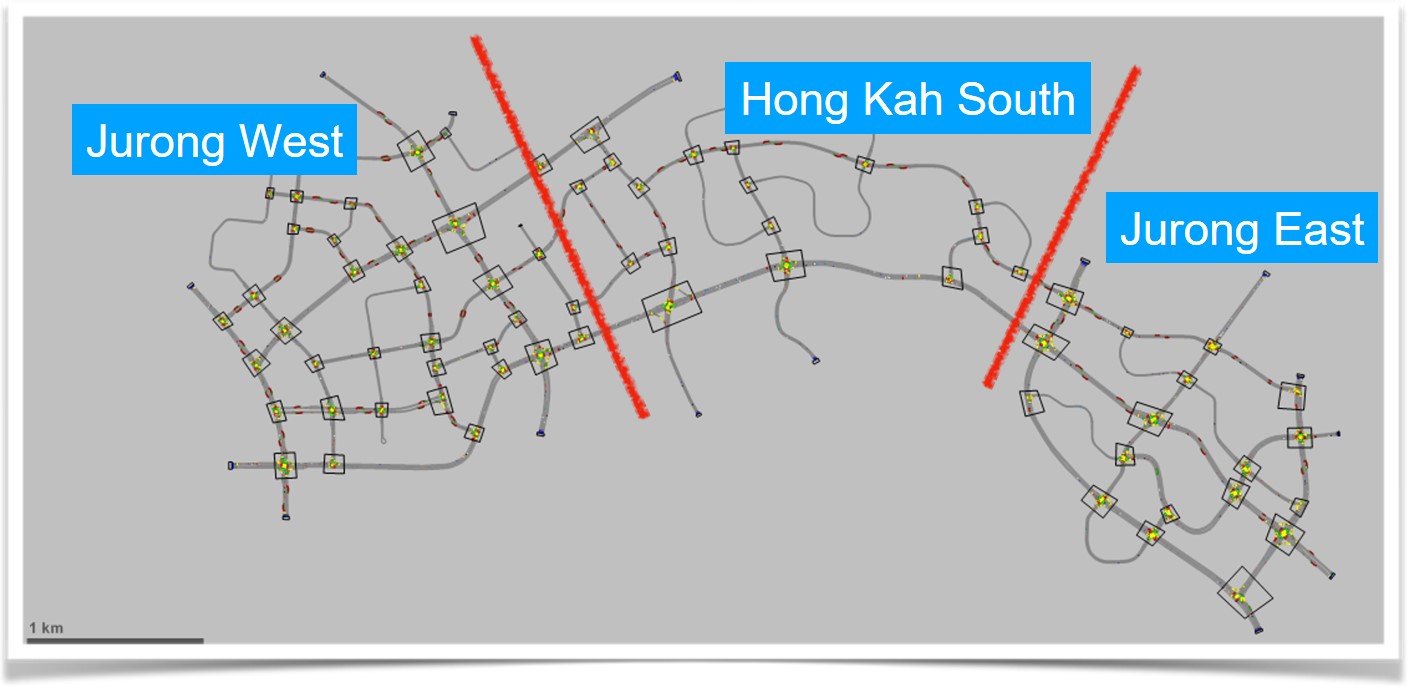}
	\caption{The traffic system settings in this case study}
	\label{fig_simusettings}
\end{figure}

\begin{figure}[!htbp]
	\centering
	\includegraphics[width=0.85\linewidth]{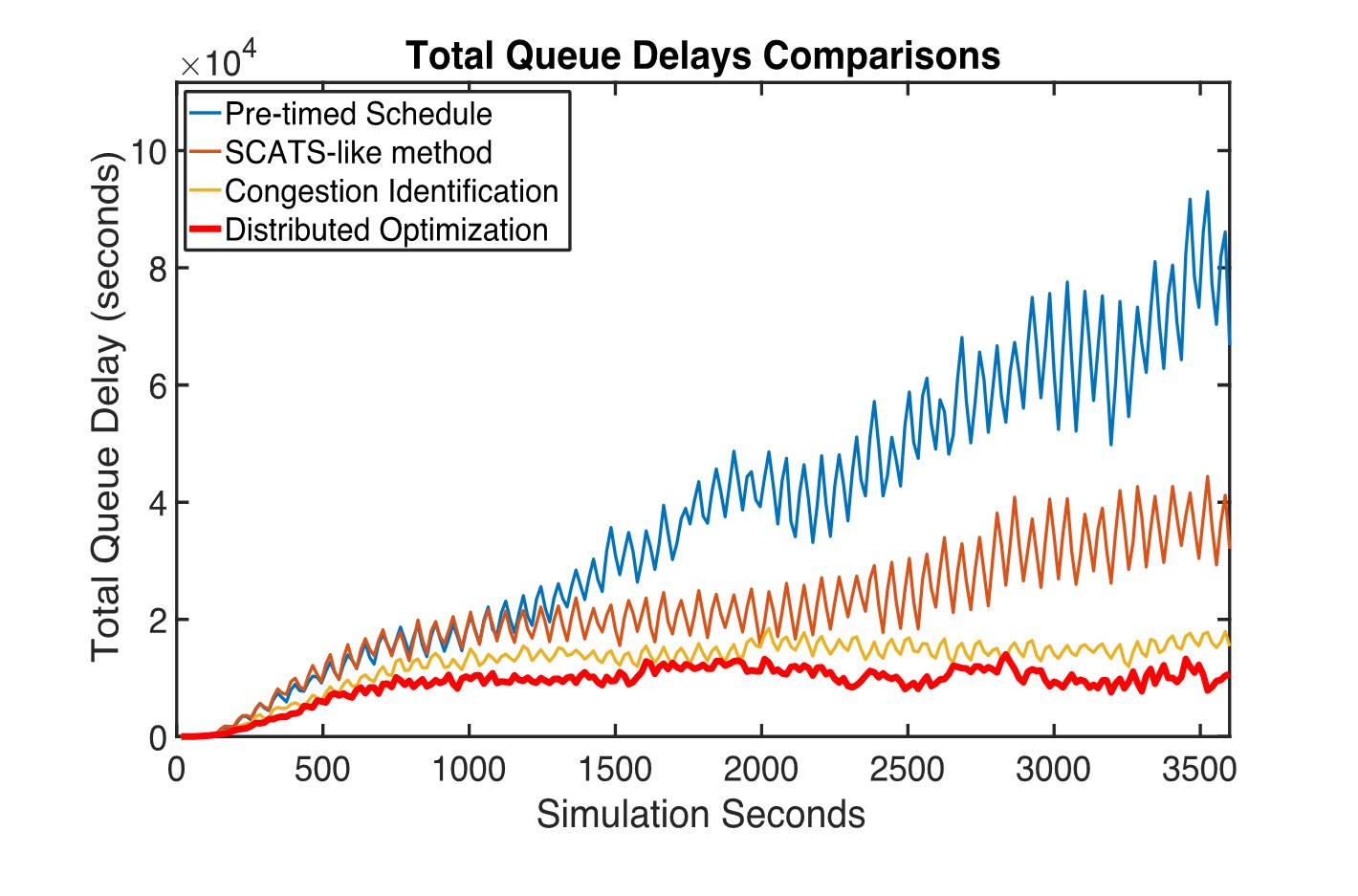}
	\caption{The simulation results on distributed traffic light control scheme based on congestion identification}
	\label{fig_simucongestion}
\end{figure}

Four curves are shown in Fig. \ref{fig_simucongestion}. We compare the proposed hybrid traffic light control strategy with two other strategies. The blue curve shows the total queue delay obtained by utilizing pre-timed scheme on the whole system, while the orange curves shows the performance obtained from a SCATS-like strategy. 
The yellow curve and red curve show the results obtained from our proposed distributed traffic light control strategy. More specifically, the red curve adopts the optimization-based traffic light control strategy, while the yellow curve shows the results obtained via the hybrid strategy. 
The simulation shows that with the proposed hybrid traffic light control strategy, the performance will be slightly degraded comparing to the performance from the optimization method. However, the computational time is reduced significantly, which allows this methodology to be implemented in a real-time manner.

\section{Related Publications}
The macroscopic traffic network modelling and the distributed computational structure have been published in \cite{Zhang2015} and \cite{zhang2017modelling}. A hybrid traffic light control strategy based on branching ratio estimation and congestion identification is presented and published in \cite{Zhang2019}. Algorithms designed based on this model have been published in \cite{Gao2017} and \cite{gao2018solving}.

\bibliography{ref}

\begin{thebibliography}{1}
\expandafter\ifx\csname url\endcsname\relax
  \def\url#1{\texttt{#1}}\fi
\expandafter\ifx\csname urlprefix\endcsname\relax\def\urlprefix{URL }\fi
\expandafter\ifx\csname href\endcsname\relax
  \def\href#1#2{#2} \def\path#1{#1}\fi

\bibitem{Zhang2015}
Y.~Zhang, R.~Su, K.~Gao, Urban road traffic light real-time scheduling, in:
  Decision and Control (CDC), 2015 IEEE 54th Annual Conference on, IEEE, 2015,
  pp. 2810--2815 (2015).

\bibitem{zhang2017modelling}
Y.~Zhang, R.~Su, C.~Sun, Y.~Zhang, Modelling and traffic signal control of a
  heterogeneous traffic network with signalized and non-signalized
  intersections, in: Control Technology and Applications (CCTA), 2017 IEEE
  Conference on, IEEE, 2017, pp. 1581--1586 (2017).

\bibitem{Zhang2019}
Y.~Zhang, Q.~Chen, R.~Su, Y.~Zhang, C.~Sun, A hybrid traffic light control
  strategy based on branching ratio estimation and congestion identification,
  in: Decision and Control (CDC), 2019 IEEE 58th Annual Conference on, IEEE,
  2019, pp. 1255--1260 (2019).

\bibitem{Gao2017}
K.~Gao, Y.~Zhang, A.~Sadollah, R.~Su, Optimizing urban traffic light scheduling
  problem using harmony search with ensemble of local search, Applied Soft
  Computing 48 (2016) 359--372 (2016).

\bibitem{gao2018solving}
K.~Gao, Y.~Zhang, R.~Su, F.~Yang, P.~N. Suganthan, M.~Zhou, Solving traffic
  signal scheduling problems in heterogeneous traffic network by using
  meta-heuristics, IEEE Transactions on Intelligent Transportation Systems
  (2018).

\end{thebibliography}

\end{document}